\documentclass[11pt,leqno]{amsart}

\usepackage{amsmath}
\usepackage{amsthm}
\usepackage{amssymb}
\usepackage{url}

\theoremstyle{definition}
\newtheorem{definition}{Definition}[section]
\newtheorem{example}[definition]{Example}
\newtheorem{remark}{Remark}
\theoremstyle{plain}
\newtheorem{thm}[definition]{Theorem}
\newtheorem{lem}[definition]{Lemma}

\begin{document}
\setcounter{page}{1} 

\vspace*{16mm}

\begin{center}
{\normalsize\bf Retarded integral inequalities of Gronwall-Bihari
type } \\[12mm]
     {\sc Rui A. C. Ferreira$^{1}$ and Delfim F. M. Torres$^{2}$}
\\ [8mm]

\begin{minipage}{123mm}
{\small {\sc Abstract.} We establish two nonlinear retarded
integral inequalities. Bounds on the solution of some retarded
equations are then obtained.}
\end{minipage}
\end{center}

 \renewcommand{\thefootnote}{}
 \footnotetext{Work supported by the \emph{Centre for Research on
Optimization and Control} and by the PhD fellowship
SFRH/BD/39816/2007 of the Portuguese Foundation for Science and
Technology.}
 \footnotetext{$^1$\,Department of Mathematics,
University of Aveiro, 3810-193 Aveiro, Portugal. E-mail:
ruiacferreira\symbol{64}ua.pt}
 \footnotetext{$^2$\,Department of Mathematics,
University of Aveiro, 3810-193 Aveiro, Portugal. E-mail:
delfim\symbol{64}ua.pt}


\vskip 12mm {\bf 1. Introduction and Preliminaries }
\medskip

\setcounter{section}{1}

In the recent paper \cite{Den} M.~Denche and H.~Khellaf study,
under some conditions on the involved functions, the following two
inequalities:
\begin{equation}
\label{in7} u(t)\leq a(t)+\int_a^t f(s)u(s)ds+\int_a^t
f(s)W\left(\int_a^s k(s,\tau)\Phi(u(\tau))d\tau\right)ds \, ,
\end{equation}
and
\begin{equation}
\label{in8} u(t)\leq a(t)+\int_a^t f(s)g(u(s))ds+\int_a^t
f(s)W\left(\int_a^s k(s,\tau)\Phi(u(\tau))d\tau\right)ds \, .
\end{equation}
Such inequalities have been then used on general time scales,
including discrete-time versions of \eqref{in7} and \eqref{in8}
(see \cite{GronTS}). In the present note we generalize both
inequalities (\ref{in7}) and (\ref{in8}) in a different direction,
by considering more general retarded inequalities, \textrm{i.e.},
by letting the upper limit of the integrals to be $C^1$
nondecreasing functions less than or equal to $t$ (\textrm{cf.}
(\ref{in1}) and (\ref{in9}) below). Moreover, our generalized
inequalities (\ref{in1}) and (\ref{in9}) are considered under less
restrictive assumptions on the involved functions, \textrm{e.g.},
in \cite{Den} the function $\Phi(\cdot)$ is assumed to be
subadditive and submultiplicative, while here we only assume
submultiplicativity. We invite the reader to compare Theorems~2.1
and 2.3 of \cite{Den} with Theorems~\ref{thm1} and \ref{thm2} of
this paper, respectively.


\bigskip
{\bf 2. Main Results}
\medskip

\setcounter{section}{2}

We start by proving a useful lemma. A similar result to
Lemma~\ref{lem1} was proved in \cite[Theorem~1.1]{Lip1} with
differentiability assumptions on the function $f(\cdot,\cdot)$.
\begin{lem} \label{lem1} Suppose that $\alpha(\cdot) \in
C^1([a,b],\mathbb{R})$ is a nondecreasing function with
$a\leq\alpha(t)\leq t$, for all $t\in[a,b]$. Assume that
$u(\cdot)$, $a(\cdot)$, $b(\cdot)\in C([a,b],\mathbb{R}_0^+)$ and
let $(t,s) \rightarrow f(t,s)\in
C([a,b]\times[a,\alpha(b)],\mathbb{R}_0^+)$ be nondecreasing in
$t$ for every $s$ fixed. If
\begin{equation*}
u(t)\leq a(t)+b(t)\int_{a}^{\alpha(t)}f(t,s)u(s)ds,
\end{equation*}
then
\begin{equation*}
u(t)\leq
a(t)+b(t)\int_{a}^{\alpha(t)}\exp\left(\int_{s}^{\alpha(t)}b(\tau)f(t,\tau)d\tau\right)f(t,s)a(s)ds\,
.
\end{equation*}
\end{lem}

\begin{proof}
The result is obvious for $t=a$. Let $t_0$ be an arbitrary number
in $(a,b]$ and define the function $z(\cdot)$ as
$$z(t)=\int_{a}^{\alpha(t)}f(t_0,s)u(s)ds,\quad t\in[a,t_0].$$
Then, $u(t)\leq a(t)+b(t)z(t)$ for all $t\in[a,t_0]$, and
$z(\cdot)$ is nondecreasing. Hence,
\begin{align*}
z'(t)&=f(t_0,\alpha(t))u(\alpha(t))\alpha'(t)\\ &\leq
f(t_0,\alpha(t))[a(\alpha(t))+b(\alpha(t))z(\alpha(t))]\alpha'(t)\\
&\leq f(t_0,\alpha(t))[a(\alpha(t))+b(\alpha(t))z(t)]\alpha'(t)\,
.
\end{align*}
The last inequality can be rearranged as
\begin{equation}
\label{eq:prv:ineq}
z'(t)-f(t_0,\alpha(t))b(\alpha(t))z(t)\alpha'(t)\leq
f(t_0,\alpha(t))a(\alpha(t))\alpha'(t).
\end{equation}
Multiplying both sides of inequality \eqref{eq:prv:ineq} by
$\exp\left(-\int_{a}^{\alpha(t)}b(s)f(t_0,s)ds\right)$, we get
\begin{multline*}
\left[z(t)\exp\left(-\int_{a}^{\alpha(t)}b(s)f(t_0,s)ds\right)\right]'\\
\leq
\exp\left(-\int_{a}^{\alpha(t)}b(s)f(t_0,s)ds\right)f(t_0,\alpha(t))a(\alpha(t))\alpha'(t)\,
.
\end{multline*}
Integrating from $a$ to $t$ and noting that $z(a)=0$, we obtain
successively that
\begin{align*}
z(t)&\leq \exp\left(\int_{a}^{\alpha(t)}b(s)f(t_0,s)ds\right)
\times \int_a^t
\exp\left(-\int_{a}^{\alpha(s)}b(\tau)f(t_0,\tau)d\tau\right) \\ &
\qquad \times f(t_0,\alpha(s))a(\alpha(s))\alpha'(s)ds\\
&=\int_a^t\exp\left(\int_{\alpha(s)}^{\alpha(t)}b(\tau)f(t_0,\tau)d\tau\right)f(t_0,\alpha(s))a(\alpha(s))\alpha'(s)ds\\
&=\int_{a}^{\alpha(t)}\exp\left(\int_{s}^{\alpha(t)}b(\tau)f(t_0,\tau)d\tau\right)f(t_0,s)a(s)ds\,
.
\end{align*}
Since $u(t)\leq a(t)+b(t)z(t)$, we have for $t=t_0$ that
$$
u(t_0)\leq
a(t_0)+b(t_0)\int_{a}^{\alpha(t_0)}\exp\left(\int_{s}^{\alpha(t_0)}
b(\tau)f(t_0,\tau)d\tau\right)f(t_0,s)a(s)ds \, .
$$
The intended conclusion follows from the arbitrariness of $t_0$.
\end{proof}

We are now in conditions to prove the following result:

\begin{thm}
\label{thm1} Suppose that $\alpha(\cdot)$, $\beta(\cdot) \in
C^1([a,b],\mathbb{R})$ are nondecreasing functions with
 $\alpha(t)$, $\beta(t) \in [a,t]$
for all $t\in[a,b]$. Assume that $u(\cdot)$, $a(\cdot)$, $b(\cdot)
\in C([a,b],\mathbb{R}_0^+)$, $(t,s) \rightarrow f(t,s)\in
C([a,b]\times[a,\alpha(b)],\mathbb{R}_0^+)$ is nondecreasing in
$t$ for every $s$ fixed, $g(\cdot,\cdot)\in
C([a,b]\times[a,\beta(b)],\mathbb{R}_0^+)$, and $(s,\tau)
\rightarrow k(s,\tau)\in
C([a,\beta(b)]\times[a,\beta(b)],\mathbb{R}_0^+)$ is nondecreasing
in $s$ for every $\tau$ fixed. Let $W(\cdot)$, $\Phi(\cdot) \in
C(\mathbb{R}_0^+,\mathbb{R}_0^+)$ be nondecreasing functions,
$\Phi(\cdot)$ submultiplicative with $\Phi(x)>0$ for $x\geq 1$.
Define
$$G(x) \doteq \int_{0}^x\frac{ds}{\Phi(1+W(s))},\quad x\geq 0,$$
$$\eta(\tau) \doteq
\max\left\{a(\tau),\int_{a}^{\beta(\tau)}g(\tau,\theta)d\theta\right\},\quad
\tau\in[a,\max\{\alpha(b),\beta(b)\}],$$
and
$$p(s) \doteq \int_{a}^{s}k(s,\tau)\Phi\left(\eta(\tau)+
b(\tau)\int_{a}^{\alpha(\tau)}\exp\left(\int_{\xi}^{\alpha(\tau)}b(\theta)f(\tau,\theta)d\theta\right)f(\tau,\xi)\eta(\xi)
d\xi \right)d\tau .$$ If for $t\in[a,b]$
\begin{equation}
\label{in1} u(t)\leq
a(t)+b(t)\int_{a}^{\alpha(t)}f(t,s)u(s)ds+\int_{a}^{\beta(t)}g(t,s)W\left(\int_{a}^{s}k(s,\tau)\Phi(u(\tau))d\tau\right)ds,
\end{equation}
then there exists $t_\ast\in(a,\beta(b)]$ such that $p(t)\in
Dom(G^{-1})$ for all $t\in[a,t_\ast]$, $G^{-1}(\cdot)$ the inverse
function of $G(\cdot)$, and
\begin{equation*}
u(t)\leq
q(t)+b(t)\int_{a}^{\alpha(t)}\exp\left(\int_{s}^{\alpha(t)}b(\tau)f(t,\tau)d\tau\right)f(t,s)q(s)ds,
\end{equation*}
where
$$q(t)=a(t)+\int_{a}^{\beta(t)}g(t,s)W\left(G^{-1}(p(s))\right)ds
.$$
\end{thm}

\begin{proof}
Let
$$
z(t)=a(t)+\int_{a}^{\beta(t)}g(t,s)W\left(\int_{a}^{s}k(s,\tau)\Phi(u(\tau))d\tau\right)ds
\, , \quad t \in [a,b] \, .
$$
Then, (\ref{in1}) can be restated as
\begin{equation}
\label{eq:h:eom} u(t)\leq
z(t)+b(t)\int_{a}^{\alpha(t)}f(t,s)u(s)ds.
\end{equation}
Applying Lemma~\ref{lem1} to \eqref{eq:h:eom}, we obtain
\begin{equation}
\label{in2} u(t)\leq
z(t)+b(t)\int_{a}^{\alpha(t)}\exp\left(\int_{s}^{\alpha(t)}b(\tau)f(t,\tau)d\tau\right)f(t,s)z(s)ds.
\end{equation}
In order to estimate $z(t)$, we define the function $v(\cdot)$ by
\begin{equation*}
v(s)=\int_{a}^{s}k(s,\tau)\Phi(u(\tau))d\tau.
\end{equation*}
We have that
$z(x)=a(x)+\int_{a}^{\beta(x)}g(x,\theta)W(v(\theta))d\theta$ and
\begin{align*}
v(s)&\leq\int_{a}^{s}k(s,\tau)\Phi\left[z(\tau) +
b(\tau)\int_{a}^{\alpha(\tau)}\exp\left(\int_{\xi}^{\alpha(\tau)}b(\theta)f(\tau,\theta)d\theta\right)f(\tau,\xi)z(\xi)
d\xi \right] d\tau\\
&\leq\int_{a}^{s}k(s,\tau)\Phi\Biggl[\eta(\tau)(1+W(v(\tau)))\\
&\quad+
b(\tau)\int_{a}^{\alpha(\tau)}\exp\left(\int_{\xi}^{\alpha(\tau)}b(\theta)f(\tau,\theta)d\theta\right)f(\tau,\xi)\eta(\xi)
d\xi(1+W(v(\tau))) \Biggr]d\tau\\
&\leq\int_{a}^{s}k(s,\tau)\Phi\Biggl[\eta(\tau)\\ &\quad +
b(\tau)\int_{a}^{\alpha(\tau)}\exp\left(\int_{\xi}^{\alpha(\tau)}b(\theta)f(\tau,\theta)d\theta\right)f(\tau,\xi)\eta(\xi)
d\xi \Biggr]\Phi(1+W(v(\tau)))d\tau\, .
\end{align*}
Let $a<t_\ast\leq \beta(b)$ be a number such that $p(t)\in Dom
(G^{-1})$ for all $t\in[a,t_\ast]$. Define $r(\cdot)$ on
$[a,s_0]$, where $a<s_0\leq t_\ast$ is an arbitrary fixed number,
by
\begin{align*}
r(s)&=\int_{a}^{s}k(s_0,\tau)\Phi\Biggl[\eta(\tau)\\ &\quad +
b(\tau)\int_{a}^{\alpha(\tau)}\exp\left(\int_{\xi}^{\alpha(\tau)}b(\theta)f(\tau,\theta)d\theta\right)f(\tau,\xi)\eta(\xi)
d\xi \Biggr]\Phi(1+W(v(\tau)))d\tau \, .
\end{align*}
Then,
\begin{align*}
r'(s)&=k(s_0,s)\Phi\Biggl[\eta(s)\\ &\quad+
b(s)\int_{a}^{\alpha(s)}\exp\left(\int_{\xi}^{\alpha(s)}b(\theta)f(s,\theta)d\theta\right)f(s,\xi)\eta(\xi)
d\xi \Biggr]\Phi(1+W(v(s)))\\ &\leq k(s_0,s)\Phi\Biggl[\eta(s)\\
&\quad +
b(s)\int_{a}^{\alpha(s)}\exp\left(\int_{\xi}^{\alpha(s)}b(\theta)f(s,\theta)d\theta\right)f(s,\xi)\eta(\xi)
d\xi \Biggr]\Phi(1+W(r(s)))\, ,
\end{align*}
that is,
\begin{multline*}
\frac{r'(s)}{\Phi(1+W(r(s)))} \leq k(s_0,s)\Phi\Biggl[\eta(s)\\ +
b(s)\int_{a}^{\alpha(s)}\exp\left(\int_{\xi}^{\alpha(s)}
b(\theta)f(s,\theta)d\theta\right)f(s,\xi)\eta(\xi) d\xi \Biggr]
\, .
\end{multline*}
Integrating both members of the last inequality from $a$ to $s$,
and having in mind that $G(r(a))=0$, we get
\begin{align*}
G(r(s))&\leq \int_{a}^{s}k(s_0,\tau)\Phi\Biggl[\eta(\tau)\\
&\quad+
b(\tau)\int_{a}^{\alpha(\tau)}\exp\left(\int_{\xi}^{\alpha(\tau)}b(\theta)f(\tau,\theta)d\theta\right)f(\tau,\xi)\eta(\xi)
d\xi \Biggr] d\tau\, .
\end{align*}
The choice of $t_\ast$ permits us to write $r(s_0)\leq
G^{-1}(p(s_0))$. Since $s_0$ is arbitrary, we conclude that (the
case $s=a$ is trivial)
\begin{equation}
\label{in3} r(s)\leq G^{-1}(p(s)),\ s\in[a,t_\ast].
\end{equation}
To complete the proof, we observe that for $a\leq s\leq t_\ast$
the inequality $\beta(\alpha(s))\leq t_\ast$ holds. Hence, we can
insert inequality (\ref{in3}) into inequality (\ref{in2}).
\end{proof}

\begin{remark}
Theorem~\ref{thm1} is new even in the particular setting studied
in \cite{Den} with $\alpha(t)=\beta(t)=t$, $b(t)=1$, and
$f(t,s)=g(t,s)=f(s)$. Indeed, one may choose in Theorem~\ref{thm1}
a submultiplicative function $\Phi(\cdot)$ that is not
subadditive, \textrm{e.g.}, $\Phi(x)=x^2$ for $x\geq 0$. This
choice of $\Phi(\cdot)$ is not a possibility in
\cite[Theorem~2.1]{Den}.
\end{remark}

To prove the forthcoming results we follow F.~M.~Dannan
\cite{Dan}, introducing the following class of functions:

\begin{definition}
\label{def} A function $g(\cdot)\in
C(\mathbb{R}_0^+,\mathbb{R}_0^+)$ is said to belong to the class
$H$ if
\begin{enumerate}
    \item $x \rightarrow g(x)$ is nondecreasing for $x\geq 0$
        and positive for $x>0$;
    \item there exists a continuous function $\Psi(\cdot)$ on
    $\mathbb{R}_0^+$ with $g(\alpha x)\leq\Psi(\alpha)g(x)$
    for $\alpha>0$, $x\geq 0$.
\end{enumerate}
\end{definition}

\begin{example}
Every continuous and nondecreasing function $g(\cdot)$ on
$\mathbb{R}_0^+$ with $g(x)>0$ for $x>0$ that is
submultiplicative, is of class $H$ with $\Psi=g$.
\end{example}

To the best of our knowledge, the following lemma is not found in
the literature. Therefore, we give a proof here.
\begin{lem}
\label{lem2} Suppose that $\alpha(\cdot)\in C^1([a,b],\mathbb{R})$
is a nondecreasing function with $a\leq\alpha(t)\leq t$ for all
$t\in[a,b]$. Assume that $u(\cdot)$, $a(\cdot)\in
C([a,b],\mathbb{R}_0^+)$ with $a(\cdot)$ a positive and
nondecreasing function, and $(t,s) \rightarrow f(t,s)\in
C([a,b]\times[a,\alpha(b)],\mathbb{R}_0^+)$ nondecreasing in $t$
for every $s$ fixed. If $g(\cdot)\in H$ and
\begin{equation}
\label{in4}  u(t)\leq a(t)+\int_{a}^{\alpha(t)}f(t,s)g(u(s))ds,
\end{equation}
then there exists a function $\Psi(\cdot)$ and a number
$t_\ast\in(a,b]$ that depends on $\Psi(\cdot)$ such that
\begin{equation}
\label{dom1}G(1)+\int_{a}^{\alpha(t)}f(t,s)\frac{\Psi(a(s))}{a(s)}ds\in
Dom(G^{-1}),\quad t\in[a,t_\ast] \, ,
\end{equation}
and
\begin{equation*}
u(t)\leq
a(t)G^{-1}\left(G(1)+\int_{a}^{\alpha(t)}f(t,s)\frac{\Psi(a(s))}{a(s)}ds\right)\,
, \quad t\in[a,t_\ast] \, ,
\end{equation*}
where
$$G(x)=\int_{x_0}^x\frac{ds}{g(s)},\quad x>0,\ x_0>0,$$
and, as usual, $G^{-1}(\cdot)$ represents the inverse function of
$G(\cdot)$.
\end{lem}

\begin{proof}
Since $a(\cdot)$ is positive and nondecreasing and $g(\cdot)\in
H$, we obtain from (\ref{in4}) that
\begin{equation*}
\frac{u(t)}{a(t)}\leq
1+\int_{a}^{\alpha(t)}\frac{f(t,s)g(u(s))}{a(s)}ds\leq
1+\int_{a}^{\alpha(t)}f(t,s)\frac{\Psi(a(s))}{a(s)}g\left(\frac{u(s)}{a(s)}\right)ds
\end{equation*}
for some function $\Psi(\cdot)$ as in the Definition~\ref{def}.
Let us now choose a number $a<t_\ast\leq b$ such that (\ref{dom1})
holds, and define function $z(\cdot)$ by
$$z(t)=1+\int_{a}^{\alpha(t)}f(t_0,s)\frac{\Psi(a(s))}{a(s)}g\left(\frac{u(s)}{a(s)}\right)ds,
\quad t\in[a,t_0],$$ where $t_0\in(a,t_\ast]$ is an arbitrary
fixed number. Then, with $x(t)=u(t)/a(t)$, we have
\begin{align*}
z'(t)&=f(t_0,\alpha(t))\frac{\Psi(a(\alpha(t)))}{a(\alpha(t))}g\left(x(\alpha(t))\right)\alpha'(t)\\
&\leq
f(t_0,\alpha(t))\frac{\Psi(a(\alpha(t)))}{a(\alpha(t))}\alpha'(t)g(z(t)),
\end{align*}
because $x(t)\leq z(t)$ and $z(t)$ is nondecreasing. Since $z(t)$
is positive, we can divide both sides of the last inequality by
$g(z(t))$ and, after integrating both sides on $[a,t]$, we get
\begin{equation*}
G(z(t))\leq
G(1)+\int_{a}^{\alpha(t)}f(t_0,s)\frac{\Psi(a(s))}{a(s)}ds\, .
\end{equation*}
Hence,
\begin{equation*}
z(t_0)\leq
G^{-1}\left(G(1)+\int_{a}^{\alpha(t_0)}f(t_0,s)\frac{\Psi(a(s))}{a(s)}ds\right)\,
.
\end{equation*}
Since $x(t_0)=u(t_0)/a(t_0)\leq z(t_0)$ and $t_0$ is arbitrary,
the result follows for all $t\in(a,t_\ast]$. The case when $t=a$
is obvious.
\end{proof}

\begin{thm}
\label{thm2} Let functions $u(\cdot)$, $f(\cdot)$, $g(\cdot)$,
$W(\cdot)$, $\Phi(\cdot)$, $\alpha(\cdot)$, $\beta(\cdot)$,
$p(\cdot)$, and $G(\cdot)$ be as in Theorem~\ref{thm1}, and
$a(\cdot)$ be as in Lemma~\ref{lem2}. If $h(\cdot)\in H$,
$$\mathcal{H}(x) \doteq
\int_{x_0}^x\frac{ds}{h(s)},\quad x>0,\ x_0>0\, ,$$ and
\begin{equation}
\label{in9} u(t)\leq a(t)+\int_{a}^{\alpha(t)}f(t,s)h(u(s))ds
+\int_{a}^{\beta(t)}g(t,s)W\left(\int_{a}^{s}k(s,\tau)\Phi(u(\tau))d\tau\right)ds
\, ,
\end{equation}
then there exists a function $\Psi(\cdot)$ and a number
$t'_\ast\in(a,\beta(b)]$ depending on $\Psi(\cdot)$ such that, for
all $t\in[a,t'_\ast]$,
$$\mathcal{H}(1)+\int_{a}^{\alpha(t)}f(t,s)\frac{\Psi(a(s))}{a(s)}ds\in
Dom(\mathcal{H}^{-1})\, ,$$
$$p(t)\in Dom(G^{-1})\, ,$$
and
$$u(t)\leq
\left[a(t)+\int_{a}^{\beta(t)}g(t,s)W\left(G^{-1}(p(s))\right)ds\right]q(t)
\, , $$
where
$$q(t) = \mathcal{H}^{-1}\left(\mathcal{H}(1)
+\int_{a}^{\alpha(t)}f(t,s)\frac{\Psi(a(s))}{a(s)}ds\right)\, .
$$
\end{thm}

\begin{proof}
Define function $z(\cdot)$ by
$$z(t)=a(t)+\int_{a}^{\beta(t)}g(t,s)W\left(\int_{a}^{s}k(s,\tau)\Phi(u(\tau))d\tau\right)ds,
\quad t\in[a,b].$$ Clearly $z(\cdot)$ is a positive and
nondecreasing function. Hence, we can apply Lemma~\ref{lem2} to
the inequality
$$u(t)\leq
z(t)+\int_{a}^{\alpha(t)}f(t,s)h(u(s))ds,$$ to obtain
\begin{equation*}
u(t)\leq z(t)
\mathcal{H}^{-1}\left(\mathcal{H}(1)+\int_{a}^{\alpha(t)}f(t,s)\frac{\Psi(a(s))}{a(s)}ds\right),
\quad t\in[a,t_\ast],
\end{equation*}
for some function $\Psi(\cdot)$ and some number $t_\ast\in(a,b]$.
An estimation of $z(t)$ can be obtained following the same
procedure as in the proof of Theorem~\ref{thm1}. After that, we
obtain
$$z(t)\leq
a(t)+\int_{a}^{\beta(t)}g(t,s)W\left(G^{-1}(p(s))\right)ds,
\quad t\in[a,t'_\ast] \, ,$$ where $G(\cdot)$ and $p(\cdot)$ are
defined as in Theorem~\ref{thm1}.
\end{proof}


\bigskip
{\bf 3. An Application}
\medskip

\setcounter{section}{3}

Let us consider the following retarded equation:
\begin{equation}
\label{eq1} u(t)=k+\int_0^{\alpha(t)}F\left(s,u(s),\int_0^s
K(\tau,u(\tau))d\tau\right)ds,\quad t\in[a,b],
\end{equation}
where $k\geq 0$, $b>0$, $\alpha(\cdot) \in C^1([a,b],\mathbb{R})$
is a nondecreasing function with $0\leq\alpha(t)\leq t$,
$u(\cdot)\in C([0,b],\mathbb{R})$, $F(\cdot)\in
C([0,b]\times\mathbb{R}\times\mathbb{R},\mathbb{R})$ and
$K(\cdot)\in C([0,b]\times\mathbb{R},\mathbb{R})$. The following
theorem gives a bound on the solution of equation (\ref{eq1}).

\begin{thm}
Assume that functions $F(\cdot,\cdot,\cdot)$ and $K(\cdot,\cdot)$
in (\ref{eq1}) satisfy
\begin{align}
|K(t,u)|&\leq k(t)\Phi(|u|)\, ,\label{in5}\\ |F(t,u,v)|&\leq
t|u|+|v| \, ,\label{in6}
\end{align}
with $k(\cdot)$ and $\Phi(\cdot)$ defined as in
Theorem~\ref{thm1}. If $u(\cdot)$ is a solution of (\ref{eq1}),
then
$$|u(t)|\leq
q(t)+t\int_{0}^{\alpha(t)}\exp\left(t(\alpha(t)-s)\right)q(s)ds,\quad
t\in[a,t_\ast],$$ for some $t_\ast\in(a,\alpha(b)]$ such that
$$p(t)\in Dom(G^{-1}),\quad t\in[a,t_\ast].$$
Here,
$$q(t)=k+\int_{0}^{\alpha(t)}G^{-1}(p(s))ds,$$
$$G(x)=\int_{0}^x\frac{ds}{\Phi(1+s)},\quad x\geq 0,$$
$$p(s)=\int_{a}^{s}k(\tau)\Phi\left[\eta(\tau)+
\tau\int_{0}^{\alpha(\tau)}\exp\left(\tau(\alpha(\tau)-\xi))\right)\eta(\xi)
d\xi \right]d\tau,$$
$$\eta(\tau)=\max\left\{k,\alpha(\tau)\right\},\quad
\tau\in[0,\alpha(b)],$$
with $G^{-1}(\cdot)$ representing the inverse function of
$G(\cdot)$.
\end{thm}

\begin{proof}
Let $u(\cdot)$ be a solution of equation (\ref{eq1}). In view of
(\ref{in5}) and (\ref{in6}), we get
\begin{equation*}
|u(t)|\leq k+\int_0^{\alpha(t)}\left(t|u(s)|+\int_0^s
k(\tau)\Phi(|u(\tau)|)d\tau\right)ds\, .
\end{equation*}
An application of Theorem~\ref{thm1} with $a(t)=k$,
$\alpha(t)=\beta(t)$, $f(t,s)=t$, $b(t)=g(t,s)=1$, and $W(u)=u$,
gives the desired conclusion:
$$|u(t)|\leq
q(t)+t\int_{0}^{\alpha(t)}\exp\left(t(\alpha(t)-s)\right)q(s)ds.$$
\end{proof}


\bibliographystyle{amsplain}

\begin{thebibliography}{99}

\bibitem{Dan} F. M. Dannan, Integral inequalities of
    Gronwall-Bellman-Bihari type and asymptotic behavior of
    certain second order nonlinear differential equations, J.
    Math. Anal. Appl. {\bf 108} (1985), no.~1, 151--164.

\bibitem{Den} M. Denche\ and\ H. Khellaf, Integral inequalities
    similar to Gronwall inequality, Electron. J. Differential
    Equations {\bf 2007}, No. 176, 14 pp. (electronic,
    \url{http://ejde.math.txstate.edu/Volumes/2007/176/abstr.html}).

\bibitem{GronTS} R. A. C. Ferreira\ and\ D. F. M. Torres,
    Generalizations of Gronwall-Bihari inequalities on time
    scales, J. Difference Equ. Appl., in press. DOI:
    10.1080/10236190802213276

\bibitem{Lip1} O. Lipovan, Integral inequalities for retarded
    Volterra equations, J. Math. Anal. Appl. {\bf 322} (2006),
    no.~1, 349--358.

\end{thebibliography}

\end{document}